\documentclass[12pt]{amsart}
\newtheorem{lemma}{Lemma}[section]
\newtheorem{theorem}[lemma]{Theorem}
\newtheorem{corollary}[lemma]{Corollary}

\usepackage{latexsym}

\begin{document}
\title{On Shapiro's Compactness Criterion for Composition Operators}
\author{John R. Akeroyd}
\date{July 9, 2010
}

\begin{abstract}
For any analytic self-map $\psi$ of $\mathbb{D} := \{z: |z| < 1\}$, J. H. Shapiro has established that the square
of the essential norm of the composition operator $C_{\psi}$ on the Hardy space $H^2$ is precisely 
$\limsup_{|w|\rightarrow 1^{-}}N_{\psi}(w)/(1 - |w|)$; where $N_{\psi}$ is the Nevanlinna counting function 
for $\psi$. In this paper we show that this quantity is equal to $\limsup_{|a|\rightarrow 1^{-}}(1 - |a|^2)||1/(1 - \bar{a}\psi)||_{H^2}^2$.
This alternative expression provides a link between the one given by Shapiro and earlier measure-theoretic notions.
Applications are given.
\end{abstract}
\maketitle
\section{Introduction and Preliminaries}
Let $\mathbb{D}$ denote the unit disk $\{z: |z| < 1 \}$ and let $\psi$ be an analytic function on $\mathbb{D}$
that maps $\mathbb{D}$ into itself; a so-called \textit{analytic self-map} of $\mathbb{D}$. Then $\psi$ has
nontangential boundary values a.e. (Lebesgue measure) on $\mathbb{T} :=\{z: |z| = 1\}$, which we also denote 
by $\psi$. Throughout this paper we let $m$ denote normalized Lebesgue measure on $\mathbb{T}$ and let $A$ 
denote normalized two-dimensional Lebesgue measure on $\mathbb{D}$. Now, by The Littlewood Subordination 
Principle (cf., \cite{Sh2}, Section 1.3), $\psi$ gives rise to a bounded composition operator $C_{\psi}$ on the 
Hardy space $H^2$ ($:= H^2(\mathbb{D})$); where $C_{\psi}(f) := f\circ\psi$.  The Nevanlinna counting function 
$N_{\psi}$ (of $\psi$) plays a central role in many results concerning $C_{\psi}$ in this and other contexts, and 
is defined on  $\mathbb{D}$ by
\[N_{\psi}(w) := -\sum_{z\in\psi^{-1}(\{w\})}\log|z|;\]
where the sum honors the multiplicity of any zero of $\psi - w$, and is zero if $\psi^{-1}(\{w\}) = \emptyset$.
Featuring among results involving $N_{\psi}$ is a theorem of J. H. Shapiro (cf., \cite{Sh1}, Theorem~2.3) that gives the square 
of the essential norm of $C_{\psi}$ on $H^2$ precisely as:
\[\limsup_{|w|\rightarrow 1^{-}}\frac{N_{\psi}(w)}{1 - |w|}.\]
Hence, $C_{\psi}$ is compact on $H^2$ if and only if 
$\lim_{|w|\rightarrow 1^{-}}\frac{N_{\psi}(w)}{1 - |w|} = 0$.  In this paper we show, by elementary methods, that 
\[\limsup_{|w|\rightarrow 1^{-}}\frac{N_{\psi}(w)}{1 - |w|} = \limsup_{|a|\rightarrow 1^{-}}\int_{\mathbb{T}}\frac{1 -  |a|^2}
{|1 - \bar{a}\psi(\zeta)|^2}\,dm(\zeta);\]
see Theorem 2.2. In the case that $C_{\psi}$ is compact on $H^2$, our alternative expression has close ties
to earlier measure-theoretic notions. We now proceed to make this connection clear. As in the introduction of \cite{Su},
let $\mu_{\psi}$ be the \textit{induced measure} of $\psi$, which is defined on Borel subsets $E$ of $\overline{\mathbb{D}}$ by
\[\mu_{\psi}(E) = m(\{\zeta\in\mathbb{T}: \psi(\zeta)\in E\}).\]
By the definition of $\mu_{\psi}$, the statement that $\lim_{|a|\rightarrow 1^{-}}\int_{\mathbb{T}}\frac{1 -  |a|^2}
{|1 - \bar{a}\psi(\zeta)|^2}\,dm(\zeta) = 0$ translates to:
\[\lim_{|a|\rightarrow 1^{-}}\int \frac{1 -  |a|^2}{|1 - \bar{a}\xi|^2}\,d\mu_{\psi}(\xi) = 0.\]
For $0 < h < 1$ and $\theta_0$ in $[0, 2\pi)$, let
\[S(h, \theta_0) = \{re^{i\theta}: 1 - h\leq r\leq 1\,\,\,\mbox{and}\,\,\, |\theta - \theta_0 | \leq h\}.\]

\begin{lemma} 
Let $\nu$ be a  finite, positive Borel measure supported in $\overline{\mathbb{D}}$. Then the following are equivalent.
\begin{itemize}
\item[$i)$] $\lim_{|a|\rightarrow 1^{-}}\int \frac{1 -  |a|^2}{|1 - \bar{a}\xi|^2}\,d\nu(\xi) = 0.$\\
\item[$ii)$] $\nu(S(h, \theta_0)) = o(h)$.
\end{itemize}
\end{lemma} 
\bigskip

\noindent Proof. By standard measure theory (cf., \cite{R}, Chapter 7), $(i)$ implies that $\nu(\mathbb{T}) = 0$, and so does $(ii)$.
Therefore, we may reduce to the case that $\nu(\mathbb{T}) = 0$. For $0 < R < 1$, let $\nu_R$ denote the restriction of $\nu$
to the annulus $\{z: R\leq |z| < 1\}$; namely, for any Borel subset $E$ of $\mathbb{D}$, $\nu_R(E) := \nu(E\cap\{z: R\leq |z| < 1\})$.
If $(i)$ holds, then, for any $\varepsilon > 0$, there exists $R$, $0 < R < 1$, such that
\[\sup_{a\in\mathbb{D}}\int \frac{1 -  |a|^2}{|1 - \bar{a}\xi|^2}\,d\nu_R(\xi) < \varepsilon.\]
Thus, by Lemma 3.3 in Chapter VI of \cite{G}, there is an absolute constant $A$ such that $\nu_R(S(h, \theta_0))\leq A\varepsilon h$,
for all $h$ and $\theta_0$. Hence, $\nu(S(h, \theta_0))\leq A\varepsilon h$ for $0 < h \leq 1 - R$ and all $\theta_0$. It follows that
$\nu(S(h, \theta_0)) = o(h)$. Conversely, suppose that $\nu(S(h, \theta_0)) = o(h)$. Then, for any $\varepsilon > 0$, there
exists $R$, $0 < R < 1$, such that $\nu_R(S(h, \theta_0))\leq \varepsilon h$, for all $h$ and $\theta_0$. Applying Lemma 3.3
in Chapter VI of \cite{G} once again, there is an absolute constant $B$ such that
\[\sup_{a\in\mathbb{D}}\int \frac{1 -  |a|^2}{|1 - \bar{a}\xi|^2}\,d\nu_R(\xi) < B\varepsilon.\]
Since $\frac{1 -  |a|^2}{|1 - \bar{a}\xi|^2}\longrightarrow 0$ uniformly on compact subsets of $\mathbb{D}$, as $|a|\rightarrow 1^{-}$,
we can now conclude that $\lim_{|a|\rightarrow 1^{-}}\int \frac{1 -  |a|^2}{|1 - \bar{a}\xi|^2}\,d\nu(\xi) = 0$.\, $\Box$
\bigskip

\noindent Now, by Lemma 1.1 and the discussion preceding it,
\[\lim_{|a|\rightarrow 1^{-}}\int_{\mathbb{T}}\frac{1 -  |a|^2}{|1 - \bar{a}\psi(\zeta)|^2}\,dm(\zeta) = 0\] 
precisely when $\mu_{\psi}(S(h, \theta_0)) = o(h)$; which is  a well-known necessary and sufficient condition 
for $C_{\psi}$ to be compact on $H^2$. This last fact is a rather straightforward exercise using the Weak Convergence 
Theorem in Section 2.4 of \cite{Sh2} along with the observation that 
\[||C_{\psi}(f)||_{H^2}^2 = \int |f|^2d\mu_{\psi};\]
once again, see the introduction of \cite{Su}, or \cite{CM} for details. 
In the next section of this paper we show that, indeed, 
$\limsup_{|w|\rightarrow 1^{-}}\frac{N_{\psi}(w)}{1 - |w|} = \limsup_{|a|\rightarrow 1^{-}}\int_{\mathbb{T}}\frac{1 -  |a|^2}
{|1 - \bar{a}\psi(\zeta)|^2}\,dm(\zeta)$ for any analytic self-map $\psi$ of $\mathbb{D}$; again, see Theorem 2.2. 
Therefore, by Shapiro's theorem, the square 
of the essential norm of $C_{\psi}$ on $H^2$ equals
\[\limsup_{|a|\rightarrow 1^{-}}\int_{\mathbb{T}}\frac{1 -  |a|^2}{|1 - \bar{a}\psi(\zeta)|^2}\,dm(\zeta)\]
(see Corollary 2.3), and $C_{\psi}$ is compact on $H^2$ precisely when 
\[\lim_{|a|\rightarrow 1^{-}}\int_{\mathbb{T}}\frac{1 -  |a|^2}{|1 - \bar{a}\psi(\zeta)|^2}\,dm(\zeta) = 0\]
(see Corollary 2.4). Applications of this are given in Corollaries 2.5 and 2.6. We close this section by recalling that B. D. MacCluer has extended the 
aforementioned measure-theoretic equivalence of the compactness of $C_{\psi}$ to the more complicated setting of the 
Hardy spaces of the ball in $\mathbb{C}^n$; cf., \cite{M}.
\pagebreak

\section{The Identity}

We begin our work here with a lemma, some notation and a few observations concerning the Nevanlinna counting 
function.

\begin{lemma}
For $0 < c < 1$,
\[(1 - c^2)\sum_{n=1}^{\infty}\frac{c^{2n-2}}{n+1}\longrightarrow 0,\]
as $c\rightarrow 1^{-}$.
\end{lemma}
\bigskip

\noindent Proof. Given $\varepsilon > 0$, there is a positive integer $N$ such that 
$\sum_{n=N}^{\infty}\frac{1}{(n+1)^2} < \varepsilon^2$, and hence:
\begin{eqnarray*}
(1 - c^2)\sum_{n=1}^{\infty}\frac{c^{2n-2}}{n+1}&\leq &
(1 - c^2)\sum_{n=1}^{N-1}\frac{c^{2n-2}}{n+1}\,\, +\,\, (1 - c^2)\left \{\sum_{n=0}^{\infty}c^{4n}\right \}^{\frac{1}{2}}
\left \{\sum_{n=N}^{\infty}\frac{1}{(n+1)^2}\right \}^{\frac{1}{2}}\\
& < & (1 - c^2)\sum_{n=1}^{N-1}\frac{c^{2n-2}}{n+1}\,\, +\,\, \left \{\sum_{n=N}^{\infty}\frac{1}{(n+1)^2}\right \}^{\frac{1}{2}} < 2\varepsilon,
\end{eqnarray*}
provided $c$ is sufficiently near $1$.\, $\Box$
\bigskip

For any point $a$ in $\mathbb{D}$, let $\varphi_a$ be the analytic 
automorphism of the unit disk given by $\varphi_a(z) = \frac{a - z}{1 - \bar{a}z}$; notice that $(\varphi_a\circ\varphi_a)(z) = z$. 
If $\psi$ is any analytic self-map of $\mathbb{D}$, then it is immediate that $N_{\psi}(\varphi_a(w)) = N_{\varphi_a\circ\psi}(w)$,
and in particular that $N_{\psi}(a) = N_{\varphi_a\circ\psi}(0)$. Parts of the proof of our main result are reminiscent of work
in Section 10.7 of \cite{Sh2}.

\begin{theorem}
Let $\psi$ be an analytic self-map of $\mathbb{D}$. Then 
\[\limsup_{|w|\rightarrow 1^{-}}\frac{N_{\psi}(w)}{1 - |w|} 
= \limsup_{|a|\rightarrow 1^{-}}\int_{\mathbb{T}}\frac{1 -  |a|^2}{|1 - \bar{a}\psi(\zeta)|^2}\,dm(\zeta).\]
\end{theorem}
\bigskip

\noindent Proof. The result clearly holds if $\psi\equiv 0$, and so we proceed under the assumption that $\psi\not\equiv 0$.
In what follows, for positive functions $f$ and $g$ of the variable $a$ in $\mathbb{D}$, we write $f(a)~\approx~g(a)$ to indicate that 
$f(a)/g(a)\longrightarrow 1$, as $|a|\rightarrow 1^{-}$. 
Then, by the Littlewood-Paley Identity and a change of variables formula 
(cf., \cite{Sh2}, pages 178 and 186, respectively), and with $c(a) := 1/|1 - \bar{a}\psi(0)|^2$,
\begin{eqnarray*}
\int_{\mathbb{T}}\frac{1 - |a|^2}{|1- \bar{a}\psi(\zeta) |^2}\,dm(\zeta)
& = & (1 - |a|^2)[c(a) + 2\int_{\mathbb{D}}\left |\frac{\bar{a}\psi '(z)}{(1- \bar{a}\psi(z))^2}\right |^2\log(\frac{1}{|z|})\,dA(z)]\\
&\approx & (1 - |a|^2)c(a) + 2\int_{\mathbb{D}}N_{\psi}(w)\frac{1 - |a|^2}{|1- \bar{a}w|^4}\,dA(w)\\
& = & (1 - |a|^2)c(a) + 2\int_{\mathbb{D}} \frac{N_{\psi}(w)}{1 - |a|^2}|\varphi_a'(w)|^2\,dA(w)\\
&\approx & (1 - |a|^2)c(a) + \int_{\mathbb{D}} \frac{N_{\psi}(w)}{1 - |a|}|\varphi_a'(w)|^2\,dA(w).\\
\end{eqnarray*} 
Since $\sup_{a\in\mathbb{D}}c(a) < \infty$, we conclude that
\[\limsup_{|a|\rightarrow 1^{-}}\int_{\mathbb{T}}\frac{1 - |a|^2}{|1- \bar{a}\psi(\zeta)|^2}\,dm(\zeta) = 
\limsup_{|a|\rightarrow 1^{-}}\int_{\mathbb{D}} \frac{N_{\psi}(w)}{1 - |a|}|\varphi_a'(w)|^2\,dA(w).\]
What remains to be shown is that 
\[\limsup_{|a|\rightarrow 1^{-}}\int_{\mathbb{D}} \frac{N_{\psi}(w)}{1 - |a|}|\varphi_a'(w)|^2\,dA(w) 
=  \limsup_{|w|\rightarrow 1^{-}}\frac{N_{\psi}(w)}{1 - |w|}.\]
To this end, first observe that
by a change of variables, the sub-averaging property (cf., \cite{Sh2}, page 190) and the notes at the beginning of this section,
\begin{eqnarray*}
\int_{\mathbb{D}} N_{\psi}(w)|\varphi_a'(w)|^2\,dA(w) &= & \int_{\mathbb{D}} N_{\psi}(\varphi_a(z))\,dA(w)\\
& = & \int_{\mathbb{D}} N_{\varphi_a\circ\psi}(z)dA(z)\\
& \geq & |\varphi_a(\psi(0))|^2N_{\varphi_a\circ\psi}(0) = |\varphi_a(\psi(0))|^2N_{\psi}(a).
\end{eqnarray*}
Thus, for $a\not = \psi(0)$,
\[\frac{N_{\psi}(a)}{1 - |a|} \leq \frac{1}{|\varphi_a(\psi(0))|^2}\int_{\mathbb{D}} \frac{N_{\psi}(w)}{1 - |a|}|\varphi_a'(w)|^2\,dA(w).\]
Since $|\varphi_a(\psi(0))|\longrightarrow 1$, as $|a|\rightarrow 1^{-}$, it now follows that
\[\limsup_{|a|\rightarrow 1^{-}}\frac{N_{\psi}(a)}{1 - |a|} \leq 
\limsup_{|a|\rightarrow 1^{-}}\int_{\mathbb{D}} \frac{N_{\psi}(w)}{1 - |a|}|\varphi_a'(w)|^2\,dA(w).\] 
For the reverse inequality, let $\beta = \limsup_{|w|\rightarrow 1^{-}}\frac{N_{\psi}(w)}{1 - |w|}$; which we know is finite
by the corollary on page~188 of \cite{Sh2}.
Then, for any $\varepsilon > 0$, there exists $R$, $0 < R < 1$, such that 
\[N_{\psi}(w) \leq (\beta + \varepsilon)(1 - |w|),\]
provided $R \leq |w| < 1$. And, by a change of variables and making use of the fact that 
$\varphi_a(z) = (a - z)\sum_{n=0}^{\infty}\bar{a}^nz^n$, we find that 
\begin{eqnarray*}
\int_{\mathbb{D}} (1 - |w|^2)|\varphi_a'(w)|^2\,dA(w) &= & 1 - \int_{\mathbb{D}} |\varphi_a(z)|^2dA(z)\\
&= & (1- |a|^2)\left [1 - (1 - |a|^2)\sum_{n=1}^{\infty}\frac{|a|^{2n - 2}}{n + 1}\right ]. 
\end{eqnarray*}
Therefore, since $|\varphi_a'(w)|^2/(1 - |a|)\longrightarrow 0$ uniformly
on compact subsets of $\mathbb{D}$, as $|a|\rightarrow~1^{-}$, and by
Lemma 2.1, we have:
\begin{eqnarray*}
\limsup_{|a|\rightarrow 1^{-}}\int_{\mathbb{D}} \frac{N_{\psi}(w)}{1 - |a|}|\varphi_a'(w)|^2\,dA(w) &\leq &
\limsup_{|a|\rightarrow 1^{-}}\int_{\mathbb{D}} \frac{(\beta + \varepsilon)(1 - |w|)}{1 - |a|}|\varphi_a'(w)|^2\,dA(w)\\
& = & \limsup_{|a|\rightarrow 1^{-}}\int_{\mathbb{D}} \frac{(\beta + \varepsilon)(1 - |w|^2)}{1 - |a|^2}|\varphi_a'(w)|^2\,dA(w)\\
& = & \limsup_{|a|\rightarrow 1^{-}}(\beta + \varepsilon)\left [1 - (1 - |a|^2)\sum_{n=1}^{\infty}\frac{|a|^{2n - 2}}{n + 1}\right ]\\
& = & \beta + \varepsilon.
\end{eqnarray*}
Now $\varepsilon > 0$ is arbitrary, and so it  follows that 
\[\limsup_{|a|\rightarrow 1^{-}}\int_{\mathbb{D}} \frac{N_{\psi}(w)}{1 - |a|}|\varphi_a'(w)|^2\,dA(w)  
\leq \limsup_{|w|\rightarrow 1^{-}}\frac{N_{\psi}(w)}{1 - |w|};\]
which completes our proof.\,$\Box$
\bigskip

Our next two results follow immediately from Theorem 2.3 in \cite{Sh1}, and Theorem 2.2 above. 
\bigskip

\begin{corollary}
Let $\psi$ be an analytic self-map of $\mathbb{D}$. Then 
\[\limsup_{|a|\rightarrow 1^{-}}\int_{\mathbb{T}}\frac{1 -  |a|^2}{|1 - \bar{a}\psi(\zeta)|^2}\,dm(\zeta)\]
is the square of the essential norm of $C_{\psi}$ on $H^2$.
\end{corollary} 
\vspace*{\fill}
\pagebreak

\begin{corollary}
Let $\psi$ be an analytic self-map of $\mathbb{D}$. Then the following are equivalent.
\begin{itemize}
\item[$i)$] $C_{\psi}$ is compact on $H^2$.
\item[$ii)$] $\lim_{|w|\rightarrow 1^{-}}\frac{N_{\psi}(w)}{1 - |w|} = 0$.
\item[$iii)$] $\lim_{|a|\rightarrow 1^{-}}\int_{\mathbb{T}}\frac{1 -  |a|^2}{|1 - \bar{a}\psi(\zeta)|^2}\,dm(\zeta) = 0$.
\end{itemize}
\end{corollary}
\bigskip

We close the paper with two applications of the equivalence of $(i)$ and $(iii)$ in Corollary~2.4. The first is well-known 
(cf., the proposition on page 32 of \cite{Sh2}), the second is less so.

\begin{corollary}
If $\psi$ is a nonconstant inner function that fixes zero, then $C_{\psi}$ is not compact on $H^2$.
\end{corollary}
\bigskip

\noindent Proof. We first observe that, for $a$ in $\mathbb{D}$ and almost all $\zeta$ in $\mathbb{T}$,
\[\frac{1}{1 - \bar{a}\psi(\zeta)} = \sum_{n=0}^{\infty}\bar{a}^n\psi^n(\zeta).\]
And since $\psi$ is a nonconstant inner function that fixes zero, $\{\psi^n\}_{n=0}^{\infty}$ 
 is an orthonomal sequence in $H^2$. Therefore,
\begin{eqnarray*}
\int_{\mathbb{T}}\frac{1 - |a|^2}{|1- \bar{a}\psi(\zeta)|^2}\,dm(\zeta) & = & (1 - |a|^2)\sum_{n=0}^{\infty}|a|^{2n}||\psi^n ||_{H^2}^2\\
&=& 1,
\end{eqnarray*}
independent of $a$ in $\mathbb{D}$.
And so, by Corollary 2.4, $C_{\psi}$ is not compact on $H^2$.\, $\Box$
\vspace*{\fill}
\pagebreak

\begin{corollary}
Let $\psi$ be an analytic self-map of $\mathbb{D}$. If $\sum_{n=0}^{\infty}||\psi^n ||_{H^2}^2$ converges, then $C_{\psi}$ is compact 
on $H^2$.
\end{corollary}
\bigskip

\noindent Proof. By our hypothesis, for any $\varepsilon > 0$, there is a positive integer $N$ such that 
$\sum_{n=N}^{\infty}||\psi^n ||_{H^2}^2 < \varepsilon^2$. Therefore, by the observation at the 
start of the proof of Corollary~2.5,

\begin{eqnarray*}
\left \{\int_{\mathbb{T}}\frac{1 - |a|^2}{|1- \bar{a}\psi(\zeta)|^2}\,dm(\zeta)\right \}^{\frac{1}{2}} 
&\leq & \sqrt{1 - |a|^2}\sum_{n=0}^{\infty}|a|^n||\psi^n ||_{H^2}\\
& \leq & \sqrt{1 - |a|^2}\left (\sum_{n=0}^{N-1}|a|^n||\psi^n ||_{H^2} + \left \{\sum_{n=0}^{\infty}|a|^{2n}\right \}^{\frac{1}{2}}
\left \{\sum_{n=N}^{\infty}||\psi^n ||_{H^2}^2\right \}^{\frac{1}{2}}\right )\\
& < & \sqrt{1 - |a|^2}\sum_{n=0}^{N-1}|a|^n||\psi^n ||_{H^2}\,\,\, + \,\,\, \varepsilon\,\,\, < \,\,\, 2\varepsilon,
\end{eqnarray*}
if $|a|$ is sufficiently near $1$. Hence, $\lim_{|a|\rightarrow 1^{-}}\int_{\mathbb{T}}\frac{1 -  |a|^2}{|1 - \bar{a}\psi(\zeta)|^2}\,dm(\zeta) = 0$. Thus,
by Corollary~2.4, $C_{\psi}$ is compact on $H^2$.\, $\Box$

\bigskip

\enlargethispage*{1000pt}
\noindent John R. Akeroyd\\
Department of Mathematics\\
University of Arkansas\\
Fayetteville, Arkansas\\
72701\\
\emph{akeroyd@uark.edu}

\end{document}